\documentclass[draft]{amsart}

\usepackage{amssymb}

\usepackage[lowtilde]{url}

\theoremstyle{plain}

\theoremstyle{definition}

\numberwithin{theorem}{section}
\numberwithin{equation}{section}

\numberwithin{enumi}{equation}

\newcommand\thmcall[1]{
  \setcounter{theorem}{\value{equation}}
  \numberwithin{equation}{theorem}    
  \numberwithin{enumi}{theorem}       
  \begin{#1}
  }

\newcommand\exitthmcall[1]{
    \end{#1}
    \numberwithin{equation}{section}  
    \numberwithin{enumi}{equation}    
    \setcounter{equation}{\value{theorem}}
  }

\newcommand\enumcall[1]{
  \stepcounter{equation}
  \begin{#1}
  }

\newcommand\exitenumcall[1]{
  \end{#1}
  }

\newcommand\bdry{\partial}
\DeclareMathOperator{\cl}{cl}
\DeclareMathOperator{\interior}{int}
\DeclareMathOperator{\LH}{LH}  

\DeclareMathOperator{\ord}{ord}
\newcommand\card{\text{card}}  

\newcommand\union{\mspace{0.5mu}\mathbin{\cup}\mspace{0.5mu}}  
\newcommand\intersect{\mathbin{\cap}}  
\newcommand\comp{\mathbin{\sim}}  

\newcommand\bigintersect{\bigcap}

\newcommand\bbz{\mathbb{Z}}
\newcommand\bbr{\mathbb{R}}
\newcommand\bbc{\mathbb{C}}
\newcommand\bbh{\mathbb{H}}

\newcommand\eitheta{e^{\imath \theta}}
\newcommand\eithetasub[1]{e^{\imath \theta_#1}}

\newcommand\dsigma{d \sigma}
\newcommand\dsimplesigma{\dtheta_1 \cdots \dtheta_n}

\newcommand\dtheta{d \theta}

\newcommand\poly[2]{P[#1,#2]}

\newcommand\abs[1]{\ensuremath{\vert #1 \vert}}
\newcommand\norm[1]{\ensuremath{\Vert #1 \Vert}}
\newcommand\inprod[2]{\ensuremath{\langle #1, #2 \rangle}}

\newcommand\plusperp{\ensuremath{\oplus_\perp}}

\newcommand\orthcomp{\ensuremath{\ominus_\perp}}

\begin{document}

\date{}

\title[An Abstract Beurling's Theorem]{An Abstract Beurling's Theorem for Several Complex Variables I}
\author{Charles W. Neville}

\address{1 Chatfield Drive, Apartment 136\\
West Hartford, CT 06110\\
USA}

\email{chip.neville at gmail.com}


\thanks{Dedicated to the memory of the late Lee Rubel.  Lee, this one's for you.}

\begin{abstract}
How to extend Beurling's theorem on the shift invariant subspaces of Hardy class $H^2$ of the unit disk to several complex variables has been an open problem at least since 1964. In this paper, we prove a generalization of Beurling's theorem to valuation Hilbert modules over valuation algebras. We shall apply our abstract Beurling's theorem to obtain a complete description of the closed invariant subspaces of $H^2$ of the polydisk.  As we shall show in a subsequent paper, further consequences of our abstract Beurling's theorem are complete extensions of Beurling's theorem to many more Hilbert spaces of analytic functions in several complex variables.
\end{abstract}

\maketitle


\setcounter{section}{0}

\section{Introduction.}\label{sec1}
When confronted with an $R$ module $M$ over a ring $R$, one of the first questions one can ask is, what are the $R$ submodules of $M$?  Often, this question is too difficult to answer, but if $M$ happens to be endowed with a topology, one can ask, what are the \emph{closed} $R$ submodules of $M$.  In a celebrated 1949 paper, Arne Beurling did just this for the case where $R$ is the algebra $\poly{\bbc}{z}$ of polynomials in the complex variable $z$, and $M$ is Hardy class $H^2$ of the unit disk \cite{Beu1}.  (See \cite{Hoff1} for complete details.)  Beurling's theorem completely characterizes the closed $\poly{\bbc}{z}$ submodules of $H^2$ of the unit disk.

As for extensions of Beurling's theorem to several complex variables, the question was raised as far back as 1964 by Henry Helson (\cite{Hel1}, p.\@ 25 ff).  In 1986, Agrawal, Clark, and Douglas observed, ``Almost everyone who has thought about this topic must have considered the corresponding problem for $H^2$ of the polydisk" but that little progress had been made, except for counter-examples and the important theorem of Ahern and Clark \cite{ACD1}.  (The theorem of Ahern and Clark characterizes the closed submodules of $H^2$ of the polydisk of finite codimension as the closures of polynomial ideals of finite codimension \cite{AC1}, cf. also the the exposition in \cite{Ru1}.)  Since Agrawal, Clark, and Douglas wrote this, significant progress was made in related areas (cf. the papers \cite{ACD1}, \cite{AB1}, \cite{Pu1}, \cite{DPSY1}, and \cite{Guo1}).  In 2007, there was a characterization of an important class of closed submodules of $H^2$ of the polydisk with deep ties to operator theory \cite{SY1}.  But a complete characterization of the closed submodules of $H^2$ of the polydisk had eluded us.. In this paper, we shall prove an abstract generalization of Beurling's classical theorem, which specializes to a concrete generalization of Beurling's theorem to $H^2$ of the polydisk.  In a second paper, we shall apply our abstract generalization of Beurlings theorem to obtain.concrete generalizations of Beurling's theorem to a number of Hilbert modules of analytic functions on domains in $\bbc^n$, and on multi-dimensional complex analytic manifolds.

\section{Some Notation.}\label{sec2}
Throughout this paper, $\bbc$ will be the complex numbers and $\bbc^n$ complex n-space.  $\bbr$ will be the real numbers and $\bbr_+$ the positive reals, $\{r \in \bbr : r \geq 0\}$.  $\bbz$ will be the integers and $\bbz_+$ the positive integers, $\{m \in \bbz : m \geq 0\}$.  $E$ will be a complex Hilbert space.  $P[z_1, \dots, z_n] = \poly{\bbc}{z_1, \dots, z_n}$ will be the algebra of polynomials with complex coefficients in the $n$ complex variables $z_1, \dots, z_n$, and $\poly{E}{z_1, \dots, z_n}$ will be the complex vector space (and $\poly{\bbc}{z_1, \dots, z_n}$ module) of polynomials with $E$-valued coefficients in the $n$ complex variables $z_1, \dots, z_n$.  We shall shamelessly identify polynomials with functions, so $\poly{E}{z_1, \dots, z_n}$ can be thought of as the space of $E$-valued polynomials in $z_1, \dots, z_n$.

We shall denote the Hilbert space norm on $E$ by $\norm{\cdot}_E$ and the inner product by $\inprod{\cdot}{\cdot}_E$.

Throughout this paper, $R$ will be a (possibly non-commutative) complex algebra (sometimes without a unit element) and $\bbh$ a left Hilbert module over $R$.  That is $\bbh$ will be both a Hilbert space and a left $R$ module, and the multiplication maps $m_r \colon h \mapsto rh$ will be continuous for each fixed $r \in R$.  We shall denote the Hilbert space norm on $\bbh$ by $\norm{\cdot}$ and the inner product by $\inprod{\cdot}{\cdot}$.

In the frequently occurring case where $\bbh$ is a Hilbert module of $E$-valued analytic functions, the norm and inner product on $E$ will be subscripted as indicated above, but they will be easy to distinguish from the corresponding quantities on $\bbh$, which will be unsubscripted.

If $X$ is a subset of a larger topological space $Y$, we shall denote the closure of $X$ by $\cl X = \cl_Y X$, the interior of $X$ by $\interior X = \interior_Y X$, and the topological boundary of $X$ by $\bdry X = \bdry_Y X$.  Finally, $\card(X)$ will denote the cardinality of the set $X$.

\section{A Preview of Coming Attractions.}\label{sec3}

The special case where the Hilbert module $\bbh$ is Hardy class $H^2$ of the polydisk and the algebra $R$ is $\poly{\bbc}{z_1, \dots, z_n}$ provides an excellent guide to our abstract Beurling's theorem.  It also illustrates why our abstract Beurling's theorem deserves to be considered a generalization of Beurling's theorem.

We need yet more standing notation.  Throughout this paper, $D^n$ will be the open $n$-dimensional unit polydisk,
\begin{equation*}
\{ z = (z_1, \dots, z_n) \in \bbc^n : \abs{z_j} < 1 \text{ for } j = 1, \dots, n \}
\end{equation*}
and $H^2(E, D^n)$ will be $E$-valued Hardy class $H^2$ of the polydisk.  Of course, $\bbc$ itself is a complex vector space (the simplest one), so $H^2(\bbc, D^n)$ is the familiar scalar valued $H^2$ space of the polydisk \cite{Ru1}.

Throughout this paper, $T^n$ will be the $n$-dimensional unit torus,
\begin{equation*}
\{ z = (z_1, \dots, z_n) \in \bbc^n : \abs{z_j} = 1 \text{ for }  j = 1, \dots, n \}.
\end{equation*}
$T^n$ is the distinguished boundary of $D^n$ and a subset of the topological boundary $\bdry D^n$ of $D^n$.  The standard parametrization of $T^n$ is $\eitheta = (\eithetasub{1}, \dots, \eithetasub{n})$, $0 \leq \theta_j \leq 2 \pi$, for $j = 1, \dots, n$.  When speaking of $T^n$, $\dsigma$ will be normalized Haar measure on $T^n$ with respect to the group of coordinate rotations,
\begin{equation*}
  z \mapsto \eitheta z = (\eithetasub{1} z_1, \dots, \eithetasub{n} z_n), \: \theta \in \bbr^n
\end{equation*}
The measure $\dsigma$ has the simple form, $\dsigma = \dsimplesigma/(2 \pi)^n$, so we shall often use the shorthand $\dsigma = \dtheta/(2 \pi)^n.$

A function $h \colon D^n \mapsto E$ is $n$-harmonic if it is continuous, and if at each point $a = (a_1, \dots, a_n) \in D^n$, each slice function
\begin{equation*}
(a_1, \dots, a_j + z_j, \dots, a_n) \mapsto h((a_1, \dots, %
a_j + z_j, \dots, a_n))
\end{equation*}
is harmonic in $z_j$ near $a_j$ (c.f.\@ \cite{Ru1}).

$H^2(E, D^n)$ is the space of all $E$-valued functions $f$ analytic on $D^n$, such that the $n$-subharmonic function $z \mapsto \norm{f(z)}_E^2$ has an $n$-harmonic majorant.  The norm on $H^2(E, D^n)$ is given by
\begin{equation*}
  \norm{f}^2 = h_f(0)
\end{equation*}
 where $h_f$ is the least $n$-harmonic majorant of the function $z \to \norm{f(z)}_E^2$.

If the Hilbert space $E$ is of finite or countable dimension, there are no measurability issues, and the well known boundary value theory for the scalar valued case (cf \cite{Ru1}) carries over directly to the $E$-valued case.  This allows us to identify $H^2(E, D^n)$ with the subspace of a.e.\@ $\dsigma$ defined functions in $L^2(E, \dsigma)$ with Fourier coefficients supported in the ``first quadrant," $\bbz_+^n$, exactly as in the scalar valued case.

However, if $E$ is of uncountable dimension, there are significant measurability issues.  Basically, these can be resolved by using the Dunford-Pettis integral and weak measurability.  ($f$ is weakly measurable if $\inprod{f}{e}_E$ is measurable for each $e \in E$.)

We need to discuss how to extend an $H^2(E, D^n)$ function $f$, only defined a.e.\@ on $T^n$, to the interior of $D^n$.  The key is to use the $n$-harmonic Poisson integral formuula
\begin{equation}
  \begin{split} \label{eqn3.1}
  f(a)   &= \int_{T^n} K(a, z) f(z) \, \dsigma(z) \\
  h_f(a) &= \int_{T^n} K(a, z) \norm{f(z)}_E^2 \, \dsigma(z)
  \end{split}
\end{equation}
for $a \in D^n$.  Here, $K(a, z)$ is the $n$-harmonic Poisson kernel, equal to the product $K_1(a_1, z_1) \cdots K_1(a_n, z_n)$ of ordinary Poisson kernels for the unit disk $D^1$, and $h_f$ is, as before, the least $n$-harmonic majorant of $\norm{f}_E^2$.  Note that $K(0,0) = 1$, so the two definitions of $H^2(E, D^n)$ are equivalent even in norm.

We can also profitably use the Szeg\H{o} integral formula for $f \in H^2(E, D^n)$:
\begin{equation}
  f(a) = \int_{T^n} S(a, z) \dsigma(z) \label{eqn3.2}
\end{equation}
for $a \in D^n$.  Here, $S(a,z)$ is the Szeg\H{o} kernel for the polydisk, equal to the product $S_1(a, z_1) \cdots S_1(a, z_n)$ of ordinary Szeg\H{o} kernels for the unit disk.

Finally, to make this absolutely simple, $H^2(E, D^n)$ is the space of square summable $E$-valued power series,
\begin{equation}
  \{ f(z) = \sum a_m z^m : \sum \norm{a_m}_E^2 < \infty \}  \label{eqn3.3}
\end{equation}
Here, and throughout this paper, we freely use standard multi-index notation:
\begin{equation*}
  a_m = a_{m_1, \dots, m_n} \quad \text{and} \quad z^m = z_1^{m_1} \cdots z_n^{m_n}
\end{equation*}
for $m = (m_1, \dots, m_n) \in \bbz_+^n$ and $z = (z_1, \dots, z_n) \in \bbc^n$.

\

In later papers, we shall consider $H^2(E, D^n(r))$, where $r$ is a real number with $0 < r < 1$, or $r = r_1, \dots r_n)$  is a poly-radius with $0 < r_j <1$ for $j =1, \dots n$.  One simply has to divide $a_j$ and $z_j$ by $r_j$ or $r$ in the above formulas to reduce the study of $H^2(E,  D^n(r))$ to that of $H^2(E, D^n))$ (\cite{Ru1} p. 3).

\

Whichever way, using formula (\ref{eqn3.1}), (\ref{eqn3.2}), or (\ref{eqn3.3}), we note several additional properties of $\bbh = H^2(E, D^n)$ and $R = \poly{\bbc}{z_1, \dots, z_n}$:
\enumcall{enumerate}
  \item $\bbh$ is a Hilbert module over the complex algebra $R$. \label{item3.4.1}
  \item Both $\bbh$ and $R$ come equipped with a $\bbz_+ \union \{\infty\}$ valued ord function which measures the order of the zero of $f \in \bbh$ or $p \in R$ at $0$. \label{item3.4.2}
  \item For $p \in R$, $\ord(p) = 0$ if $p$ is a unit of $R$. \label{item3.4.3}
  \item For $f \in \bbh$ or $f \in R$, $\ord(f)  = \infty$ if and only if $f$ is identically $0$. \label{item3.4.4}
  \item The  $\ord$ function almost respects multiplication: $\ord(p f) \geq \ord(p)$\\ $+ \ord(f)$, and $\ord(p q) \geq \ord(p) + \ord(q)$, for $f \in \bbh$, $p \in R$, and $q \in R$.  (Of course, in this case, inequality may be replaced by equality, but we want to state this in as general a form as possible.) \label{item3.4.5}
  \item The $\ord$ function respects scalar multiplication: $\ord(\lambda f) = \ord(f)$ for $\lambda \in \bbc$, $\lambda \neq 0$, and $f \in \bbh$ or $f \in R$. \label{item3.4.6}
  \item The $\ord$ function (almost) respects addition: $\ord(f + g) \geq \\
   \min( \ord(f), \ord(g))$ for $f \text{ and } g \in \bbh$ or $f \text{ and } g \in R$. \label{item3.4.7}
  \item The $\ord$ function is upper semi-continuous on $\bbh$. Thus the subspaces $\bbh_m = \{ f \in \bbh : \ord(f) \geq m \}$ are closed for $m = 1, 2, \dots$ because $\bbh_m$ is the inverse image of the set $R_m  = \{ r \in R : \ord(r) \geq m\}$. \label{item3.4.8}
\exitenumcall{enumerate}

There are, in fact, several possible $\ord$ functions.  The one we shall use is

\thmcall{definition}\label{def3.5}
$\ord(f) = $ the degree of the lowest degree term in the power series expansion of $f$ for $f \in \bbh$, and similarly for $R$.
\exitthmcall{definition}

We still need to define the degree of a monomial:  The \emph{degree} of $z^m = z_1^{m_1} \cdots z_n^{m_n}$ is $\abs{m}$, where here and throughout this paper, $\abs{m} = \abs{m_1} + \dots + \abs{m_n}$.  This, of course, is the standard definition.

Another equivalent definition of the $\ord$ function is

\thmcall{definition}\label{def3.6}
$\ord(f) = $ the lowest degree term in the expansion of $f$ in homogeneous polynomials for $f \in \bbh$, and similarly for $p \in R$.
\exitthmcall{definition}

Of course, the degree of a homogeneous polynomial is the common degree of its monomial terms.  (See \cite{Ru2} for homogeneous expansions of analytic functions.)

\

In the next section, we shall call a left Hilbert module $\bbh$ over a (possibly non-commutative) complex algebra $R$ a \emph{valuation Hilbert module}, and $R$ a \emph{valuation algebra}, if both are equipped with an $\ord$ function satisfying properties (\ref{item3.4.1}) through (\ref{item3.4.8}).  Our abstract Beurling's theorem will characterize the closed left $R$ submodules of of an arbitrary valuation Hilbert module.  (Truth in advertising: Actually, the closed $R_1$ submodules, where $R_1 = \{ r \in R : \ord(r) \geq 1 \}$.  But if $R$ has a unit element $1_R$ and $R = \bbc \cdot 1_R  \,+\,  R_1$, these are the same as the closed left $R$ submodules.  See definition \ref{def4.10} and proposition \ref{prop4.11}.)  As hinted above, homogeneous subspaces, similar to subspaces of of homogeneous polynomials all of the same degree or $\{0\}$, will play a prominent role.

\

To gain an inkling of what our actual abstract Beurling's theorem will look like in the polydisk case, recall the wandering subspace approach to the classical $1$-variable Hilbert space valued Beurling's theorem due to Halmos \cite{Hal1}.  (See \cite{Hoff1} for a complete exposition.)  This theorem characterizes the closed $\poly{\bbc}{z}$ submodules (classically called closed shift invariant subspaces) of $H^2(E, D^1)$ as the subspaces of the form
\begin{equation}
  V = I \cdot H^2(W, D^1) \label{eqn3.7}
\end{equation}
where $W$ is a complex Hilbert space and $I$ is an operator valued inner function.  This means
\begin{equation}
  I(z) \colon W \mapsto E \label{eqn3.8}
\end{equation}
is a (usually partial) isometry a.e.\@ $\dtheta$ on $T^1$, the unit circle.  Given a closed $\poly{\bbc}{z}$ submodule $V$, the operator valued inner function is constructed from the \emph{wandering subspace}
\begin{equation}
  W = V \orthcomp z V \label{eqn3.9} 
\end{equation}
Here and throughout this paper, $\orthcomp$ will denote the orthogonal complement.  Then an equivalent statement of the $1$-variable Beurling's theorem is
\begin{equation}
  V = W \plusperp z W \plusperp z^2 W \plusperp \cdots \label{eqn3.10}
\end{equation}
where here and throughout this paper, $\plusperp$ denotes the orthogonal direct sum.

\

Our abstract Beurling's theorem uses a similar idea.  Specialized to the case of $H^2(E, D^n)$, our theorem may be briefly stated as follows:

Let $V$ be a closed subspace of $H^2(E, D^n)$.  In place of the single wandering subspace $W$ in equations (\ref{eqn3.9}) and (\ref{eqn3.10}), we consider the decreasing sequence of closed subspaces
\begin{equation}
\begin{split} \label{eqn3.11}
  V_0 = V, \quad V_1 &= \{ f \in V : \ord(f) \geq 1 \},  \\
                 V_2 &= \{ f \in V : \ord(f) \geq 2 \}, \: \dots
\end{split}
\end{equation}
and the sequence of \emph{near homogeneous} subspaces,
\begin{equation}
  W_0 = V_0 \orthcomp V_1, \quad W_1 = V_1 \orthcomp V_2, \: \dots \label{eqn3.12}
\end{equation}
(A subspace is \emph{near homogeneous} if the $\ord$ function is constant on all non-zero elements of the subspace.)  Then (as we shall prove later although it is fairly obvious)
\begin{equation}
  V = W_0 \plusperp W_1 \plusperp W_2 \plusperp \cdots \label{eqn3.13}
\end{equation}
We call this orthogonal decomposition the \emph{near homogeneous decomposition} of $V$.

Any closed subspace has a near homogeneous decomposition, but closed left submodules also have a \emph{near inner} decomposition.  In the polydisk case, a closed subspace $W$ of $H^2(E, D^n)$ is \emph{near inner} if $W \perp z^m W$ for all multi-indices $m \in \bbz_+^n$ with $\abs{m} = \abs{m_1} + \cdots + \abs{m_n} \geq 1$.  We can define a \emph{near inner decomposition} of a closed subspace $V$ to be an orthogonal decomposition
\begin{equation}
  V = W^0 \plusperp W^1 \plusperp W^2 \plusperp \cdots \label{eqn3.14}
\end{equation}
of near inner subspaces such that
\begin{equation}
  z^k W^l \perp W^m  \label{eqn3.15}
\end{equation}
whenever $\abs{k} + l > m$.  (The actual definition is slightly more complicated -- see definition \ref{def6.3}.)  We have superscripted the subspaces $W^0, W^1, W^2, \ldots$ to distinguish them from the near homogeneous decomposition of $V$.  However, there is a close connection between near homogeneous decompositions and near inner decompositions, and this connection constitutes the core of our abstract Beurling's theorem:

\thmcall{theorem} \label{thm3.16}
Let $V$ be a closed subspace of $H^2(E, D^n)$.  Then $V$ is a closed $\poly{\bbc}{z_1, \dots, z_n}$ submodule of $V$ if and only if the near homogeneous decomposition of $V$ is near inner and has the full projection property.
\exitthmcall{theorem}

The full projection property, which we shall not go into until section \ref{sec7}, is designed to eliminate situations such as
\begin{equation}
  V = \LH\{I\} \plusperp \LH\{z^2 I\} \plusperp \LH\{z^4 I\} \plusperp \cdots \label{eqn3.17}
\end{equation}
where $I$ is an inner function in $H^2(\bbc, D^1)$.  (Here and throughout this paper, $\LH$ denotes the linear hull of a set.)  In this case, the orthogonal decomposition (\ref{eqn3.17}) is the near homogeneous decomposition of $V$, and it is near inner, but it does not have the full projection property.  This is a good thing because $V$ is not a $\poly{\bbc}{z}$ submodule of $H^2(\bbc, D^1)$.

\

In the classical $1$-variable case, our theorem \ref{thm3.16} says the the closed $\poly{\bbc}{z}$ submodules of $H^2(\bbc, D^1)$ are exactly the closed subspaces of the form (\ref{eqn3.10}), although this takes some proving.  In fact, the decomposition (\ref{eqn3.10}) turns out to be the near homogeneous decomposition of $V$.  Furthermore, it is near inner and has the full projection property, as we shall show in a second paper.  Thus our abstract Beurling's theorem really is a generalization of the classical $1$-variable Beurling's theorem in the case where the Hilbert space of values equals $\bbc$.

\

We shall close with a brief outline of things to come.  In section \ref{sec4}, we shall carefully develop the theory of valuation Hilbert modules over valuation algebras.  In section \ref{sec5}, we shall discuss near homogeneous subspaces and decompositions.  In section \ref{sec6}, we shall discuss near inner subspaces and decompositions, and give the full, correct definition of a near inner decomposition.  In section \ref{sec7}, we shall discuss the full projection property.  In section \ref{sec8}, we shall fully state and prove our abstract Beurling's theorem characterizing the closed submodules of valuation Hilbert modules over valuation algebras.  And in section \ref{sec9}, we shall apply our abstract Beurling's theorem to completely characterize the closed invariant subspaces of $H^2(E, D^n)$.

In our second paper mentioned above, we shall discuss analytic algebras and analytic Hilbert modules, which provide a general setting for examples where our abstract Beurling's theorem applies.  Our second paper paper will then present concrete examples of the application of our abstract Beurling's theorem to several complex variables, including to $H^2$ spaces and Bergman $A^2$ spaces. And at the end of our second paper, we shall list of a number of open problems.

\section{Valuation Algebras and Hilbert Modules.}\label{sec4}
The notion of a valuation on a field will probably be familiar to most readers.  An extension of one form of this notion to algebras and Hilbert modules will provide the setting for our abstract Beurling's theorem.  As you may have already guessed from section \ref{sec3}, the motivating example for a \emph{valuation algebra} and \emph{valuation Hilbert module} is the algebra $\poly{\bbc}{z_1, \dots, z_n}$ and Hilbert module $H^2(E, D^n)$, with the $\ord$ function introduced there as the valuation.

There are several notions of a valuation in abstract algebra.  The one we shall borrow (and amend slightly) is often called an \emph{exponential valuation} \cite{Hus1}.  Perhaps we should call our notion of a valuation an`\emph{inverse exponential valuation}, but the phrase is way too cumbersome, so we shall simply use the term \emph{analytic valuation} to avoid confusion and distinguish it from the usual notion of a valuation in abstract algebra.  For concision, we shall usually drop the word {analytic} and smply call our notion a \emph{valuation}.

We state here and now that ALL algebras and Hilbert spaces will be \emph{complex}.

\thmcall{definition}\label{def4.1}
Let $R$ be a (possibly non-commutative) algebra.  An \emph{analytic valuation} on $R$ is a function $\ord \colon \! R \mapsto \bbz_+ \union \{\infty\}$ such that for all $r$ and $s \in R$,
\enumcall{enumerate}
  \item $\ord(r) = 0$ if $r$ is a left or right unit of $R$. \label{item4.1.1}
  \item $\ord(r) = \infty$ if and only if $r = 0$. \label{item4.1.2}
  \item $\ord(r s) \geq \ord(r) + \ord(s)$.  \label{item4.1.3}
  \item $\ord(\lambda r) = \ord(r)$ for $\lambda \in \bbc, \lambda \neq 0$ \label{item4.1.4}
  \item $\ord(r + s) \geq \min(\ord(r), \ord(s))$  \label{item4.1.5}
\exitenumcall{enumerate}
\exitthmcall{definition}
Of course, condition (\ref{item4.1.1}) is satisfied vacuously if $R$ does not have a two sided identity element.

Here and throughout, we follow the usual conventions with respect to $\infty$: $m < \infty$ for all $m \in \bbz$, and for such $m$,  $m \cdot \infty = \infty \cdot m = \infty$ if $m \neq 0$.  Further, $\infty \cdot 0 = 0 \cdot \infty = \infty$ and $\infty \cdot \infty = \infty$.

\thmcall{definition}\label{def4.2}
A \emph{valuation algebra} is an ordered pair $(R, \ord)$, where $R$ is a (possibly non-commutative) algebra and $\ord$ is an analytic valuation on $R$.
\exitthmcall{definition}

\thmcall{definition}\label{def4.3}
Let $(R, \ord_R)$ be a valuation algebra, and let $\bbh$ be a Hilbert space which is a left $R$ Hilbert module.  A \emph{Hilbert module valuation} on $\bbh$ (with respect to $(R, \ord_R)$) is a function $\ord \colon \bbh \mapsto \bbz_+ \union \{\infty\}$ such that for all $h$, $h_1$ and $h_2 \in \bbh$, and for all $r \in R$,
\enumcall{enumerate}
  \item $\ord(h) = \infty$ if and only if $h = 0$. \label{item4.3.1}
  \item $\ord(r h) \geq \ord_R(r) + \ord(h)$  \label{item4.3.2}
  \item $\ord(\lambda h) = \ord(h)$ for $\lambda \in \bbc, \lambda \neq 0$. \label{item4.3.3}
  \item $\ord(h_1 + h_2) \geq \min(\ord(h_1), \ord(h_2))$. \label{item4.3.4}
  \item The $\ord$ function is upper semi-continuous on $\bbh$. \label{item4.3.5}
\exitenumcall{enumerate}
\exitthmcall{definition}

Note that each set $\bbh_m$ is a subspace by properties (\ref{item4.3.3}) and (\ref{item4.3.4}).  The subspaces $\bbh_m = \{ h \in \bbh : \ord(h) \geq m \}$ are closed for $m = 0, 1, 2, \, \dots$ because $\ord$ is upper semicontinuous. The subspace $\bbh_0$ is automatically a closed subspace because $\bbh_0 = \bbh$. Fiinally, as part of the definition of an $R$ Hilbert module, the multiplication map $m_r \colon h \mapsto r h$ is continuous for each $r \in R$.

We need to review some of the properties of upper semicontinuous functions (\cite{Wiki1},  \cite{H1}): Here is the easiest way to see that the prototype of the $\ord$ function is upper semicontinuous at $h$, but that we cannot define $\ord$ to be a continuous function: Recall that the one of the conditions that $\ord$ is upper semicontinuous at $h \in \bbh$ is 
    \begin{equation*}
        \limsup \ord(h_k) \leq \ord(h) \text{ if } h_k \text{ is a sequence in } \bbh \text{ converging to } h.
    \end{equation*}
Now suppose that $\bbh = H^2$, and recall that for $h \in H^2$, $\ord( h )$ is the order of the zero of $h$ at $0$. The upper semicontinuity condition is clearly satisfied by $\ord$ in $H^2$.

Now if $h(0) \neq 0$, $\ord( h_k )$ converges to $\ord( h ) = 0$, so $\ord$ is continuous at $h$. But if $h$ has a zero at $0$, then $\ord( h ) > 0$. Choose $h_k = h + 1/k$ and observe that $h_k$ converges to $h$, but $\ord( h_k ) = 0$ does not converge to $\ord( h )$. Thus $\ord$ is not continuous at $h$.

\thmcall{definition}\label{def4.4}
Let $(R, \ord_R)$ be a valuation algebra.  A \emph{valuation Hilbert module} over $(R, \ord_R)$ is an ordered pair $(\bbh, \ord)$, where $\bbh$ is a left Hilbert module over $R$ and $\ord$ is a Hilbert module valuation on $\bbh$ with respect ot $(R, \ord_R)$.
\exitthmcall{definition}

Where it will cause no confusion, we shall denote both the algebra valuation on $R$ and the Hilbert module valuation on $\bbh$ by the same symbol, $\ord$.  We shall also make the gloss, whenever convenient, of denoting the valuation algebra $(R, \ord)$ by $R$ alone, and the valuation Hilbert module $(\bbh, \ord)$ by $\bbh$ alone.  Finally, we note that the concepts of a valuation algebra and a valuation Hilbert module apply equally well to \emph{real} algebras and \emph{real} Hilbert modules, but we shall not consider these here.  Our focus will be exclusively on \emph{complex} valuation algebras and Hilbert modules.

The inequalities in properties (\ref{item4.1.3}) and (\ref{item4.3.2}) will doubtless raise the eyebrows of some readers, as will the fact that the algebra $R$ might be non-commutative, might have divisors of zero, and might not have an identity element.  Normally in abstract algebra, a valuation algebra must be an integral domain (ie.\@ be commutative and have no divisors of zero) with an identity element, and there must be equality rather than inequality in property (\ref{item4.1.3}).  There are examples of valuation Hilbert modules over valuation algebras of bounded operator valued functions which are non-commutative and have divisors of zero. In these examples, there will sometimes be inequality in properties (\ref{item4.1.3}) and (\ref{item4.3.2}). However, we shall not consider these in this or our subsequent paper.
 
In the most common examples of valuation algebras and valuation Hilbert modules, equality does hold.  We need a name for this situation:

\thmcall{definition}\label{def4.5}
An algebra valuation $\ord$ on an algebra $R$ is \emph{strict} if equality holds in property (\ref{item4.1.3}), that is if
\begin{equation}
  \ord(r s) = \ord(r) + \ord(s) \label{eqn4.5.1}
\end{equation}
for all $r$ and $s \in R$.  A \emph{strict} valuation algebra is a valuation algebra for which the valuation $\ord_R$ is strict.

A Hilbert module valuation $\ord$ on a left Hilbert module $\bbh$ is \emph{strict} if equality holds in property (\ref{item4.3.2}), that is if
\begin{equation}
  \ord(r h) = \ord(r) + \ord(h) \label{eqn4.5.2}
\end{equation}
for all $r \in R$ and $h \in \bbh$.  A \emph{strict} valuation Hilbert module is a valuation Hilbert module for which the valuation $\ord$ is strict.
\exitthmcall{definition}

If an algebra $R$ or a left Hilbert module over $R$ admits a strict valuation, then the algebra must be quite well behaved.

\thmcall{proposition}\label{prop4.6}
Let $R$ be a strict valuation algebra.  Then $R$ has no left divisors of zero and no right divisors of zero.  If $R$ is commutative, then it is an integral domain.
\exitthmcall{proposition}

\begin{proof}
Suppose $R$ had divisors of zero.  Let $r$ be a left divisor of zero and $s$ be a right divisor of zero, such that $r \neq 0, s \neq 0$, but $r s = 0$.  Then $\ord(r s) = \infty > \ord(r) + \ord(s)$, which contradicts the strictness of $R$.
\end{proof}

Note that the second assertion of the proposition, that $R$ is an integral domain, follows immediately from the definition of an integral domain as a commutative algebra with no divisors of zero.

\thmcall{proposition}\label{prop4.7}
Let $\bbh$ be a strict valuation Hilbert module over a valuation algebra $R$.  Then $R$ is a strict valuation algebra and has no left divisors of zero, and no right divisors of zero.  Further, if $R$ is commutative, then $R$ is an integral domain.
\exitthmcall{proposition}

\begin{proof}
That $R$ has no left divisors of zero and no right divisors of zero follows immediately from proposition \ref{prop4.6} once we prove that $R$ is strict.  To that end, let $r$ and $s$ be non-zero elements of $R$, and let $h$ be a non-zero element of $\bbh$.  Then
\begin{equation}
\begin{split} \label{eqn4.8}
  \ord((r s) h) &= ord(r s) + ord(h) \\
  \ord((r s) h) &= \ord(r (s h)) = ord(r) + ord(s h) \\
                &= \ord(r) + \ord(s) + \ord(h)
\end{split}
\end{equation}
Thus $\ord((r s) h) \neq \infty$, and $\ord(h) \text{\ certainly} \neq \infty$, so we can set
\begin{equation*}
  ord(r s) + ord(h) = \ord(r) + \ord(s) + \ord(h)
\end{equation*}
and cancel out $\ord(h)$ to obtain $\ord(r s) = \ord(r) + \ord(s)$.  Hence $R$ is strict.
\end{proof}

The subalgebra $R_1 = \{ r \in R : \ord(r) \geq 1 \}$ will be particularly important in this paper.  In fact, the natural form of our abstract Beurling's theorem will concern $R_1$ invariant subspaces.

\thmcall{definition}\label{def4.9}
Let $R$ be a valuation algebra and let $\bbh$ be a valuation Hilbert module over $R$.  Then a closed subspace $V$ of $\bbh$ is $R_1$ \emph{invariant} if $R_1 V \subseteq V$.
\exitthmcall{definition}

A natural question is, when is an $R_1$ invariant subspace $V$ of $\bbh$ a left submodule, that is when is $R V \subseteq V$?  A useful condition which guarantees this is given by the following definition and proposition:

\thmcall{definition}\label{def4.10}
Let $R$ be a valuation algebra with unit element (ie. identity element) $1_R$.  Then $R$ is a \emph{unital} valuation algebra if
\begin{equation}
  R = \{ \lambda 1_R : \lambda \in \bbc \} + R_1  \label{eqn4.10.1}
\end{equation}
\exitthmcall{definition}

Note that the $+$ in equation \ref{eqn4.10.1} represents simple addition.  It does not represent a direct sum.

\thmcall{proposition}\label{prop4.11}
Let $R$ be a valuation algebra and let $\bbh$ be a left Hilbert module over $R$.  Then every $R_1$ invariant subspace of $\bbh$ is a left $R$ submodule of $\bbh$ if either of the following conditions hold:
\enumcall{enumerate}
\item  $R = R_1$, or \label{item4.11.1}
\item  $R$ is unital. \label{item4.11.2}
\exitenumcall{enumerate}
\exitthmcall{proposition}

\begin{proof}
Clear.
\end{proof}

A simple example of a valuation algebra satisfying condition \ref{item4.11.1} is
\begin{equation*}
P_1[\bbc, z_1, \dots, z_n] = \{ f \in \poly{\bbc}{z_1, \dots, z_n} : f(0) = 0 \}
\end{equation*}
Two iconic examples of valuation algebras satisfying condition \ref{item4.11.2} are the algebra of polynomials $\poly{\bbc}{z_1, \dots, z_n}$ and the algebra $H^\infty(\bbc, (M, p))$ of bounded complex valued analytic functions on a complex manifold $M$ with basepoint $p$.  In the first case, $\ord(f) =$ the order of the zero of $f$ at $0$, and in the second $\ord(f) =$ the order of the zero of $f$ at $p$.

\section{Near Homogeneous Decompositions.}\label{sec5}
The additional structure provided by the valuation on a valuation Hilbert module allows us to decompose its elements into parts analogous to the decomposition of a power series into a sum of homogeneous polynomials.  We shall call these decompositions \emph{near homogeneous decompositions}.

It will be convenient to introduce some more standing notation.  Unless otherwise stated, $R$ will be a valuation algebra, $\bbh$ will be a valuation Hilbert module over $R$, and $V$ will be a closed subspace of $\bbh$.  To avoid trivial issues, unless otherwise stated, $V \neq \{ 0 \}$.
  Remember from our standing notation in section \ref{sec2} and definition \ref{def4.4} that if $R$ is non-commutative, then $\bbh$ is a \emph{left} (as opposed to a \emph{right}) Hilbert module over $R$.

\thmcall{definition}\label{def5.1}
The valuation subspace series for $V$ is the decreasing sequence of closed subspaces
\begin{equation}
  V = V_0 \supseteq V_1 \supseteq V_2 \supseteq \cdots \label{eqn5.1.1}
\end{equation}
where
\begin{equation}
  V_k = \{ h \in V : \ord(h) \geq k \} \label{eqn5.1.2}
\end{equation}
The valuation ideal series for $R$ is the decreasing sequence of two sided ideals
\begin{equation}
  R = R_0 \supseteq R_1 \supseteq R_2 \supseteq \cdots \label{eqn5.1.3}
\end{equation}
where
\begin{equation}
  R_k = \{ r \in R : \ord(r) \geq k \} \label{eqn5.1.4}
\end{equation}
\exitthmcall{definition}

We have already met the valuation subspace series for $\bbh$; it is the decereasing sequence of closed subspaces $\bbh = \bbh_0 \supseteq \bbh_1 \supseteq \bbh_2 \supseteq \cdots$ following property (\ref{item4.3.5}).  Note that each of the sets $V_k$ is a subspace by properties (\ref{item4.3.3}) and (\ref{item4.3.4}).  Because $V_k = \bbh_k \intersect V$, each of the sets $V_k$ in the valuation subspace series (\ref{eqn5.1.1}) for $V$ is a closed subspace of $V$.

That each of the sets $R_k$ in the valuation ideal series (\ref{eqn5.1.3}) is in fact a two sided ideal follows from properties (\ref{item4.1.3}), (\ref{item4.1.4}), and (\ref{item4.1.5}).  That
\begin{equation}
  \bigintersect R_k = \{ 0 \}, \quad \bigintersect V_k = \{ 0 \}, \quad \bigintersect \bbh_k = \{ 0 \} \label{eqn5.2}
\end{equation}
all follow from properties (\ref{item4.1.2}) and (\ref{item4.3.1}) To see this, suppose, that $h \in \bigintersect V_k$. Then for each $k$, $h \in V_k$, so $\ord(h) \geq k$ for each $k < \infty$. Thus $\ord(h) = \infty$, so $h = 0$ by property (\ref{item4.3.1}).

The valuation ideal series for $R$ is entirely analogous to the valuation subspace series for $V$, except that we can use the Hilbert space structure of $V$ to form orthogonal complements:

\thmcall{definition}\label{def5.3}
The \emph{near homogeneous decomposition} of $V$ is the orthogonal direct sum
\begin{equation}
  V = W_0 \plusperp W_1 \plusperp W_2 \plusperp \, \cdots \label{eqn5.3.1}
\end{equation}
where
\begin{equation}
  W_0 = V_0 \orthcomp V_1, \; W_1 = V_1 \orthcomp V_2, \; W_2 = V_2 \orthcomp V_3, \, \dots \label{eqn5.3.2}
\end{equation}
and $V = V_0 \supseteq V_1 \supseteq V_2 \supseteq \, \cdots$ is the valuation subspace series for $V$.  The orthogonal subspaces $W_0$, $W_1$, $W_2, \, \dots$ are called the \emph{near homogeneous components} of $V$.
\exitthmcall{definition}

We need an all too obvious comment here.  For algebraists, a direct sum like (\ref{eqn5.3.1}) is the subspace of all \emph{finite} sums of elements from the subspaces $W_k$, $k \in \bbz_+$.  But for analysts like us, the orthogonal direct sum (\ref{eqn5.3.1}) is the closed subspace of all norm convergent sums, whether \emph{finite} or \emph{infinite}, of elements from the subspaces $W_k$, $k \in \bbz_+$.

Note that the near homogeneous components $W_0, W_1, W_2, \, \dots$ are closed because they are orthogonal complements of sets in closed subspaces.  Further, they are orthogonal by construction.  However, that the orthogonal direct sum (\ref{eqn5.3.1}) equals all of $V$ requires proof:

\thmcall{proposition}\label{prop5.4}
The orthogonal direct sum
\begin{equation*}
  W = W_0 \plusperp W_1 \plusperp W_2 \plusperp \, \cdots
\end{equation*}
equals all of $V$.
\exitthmcall{proposition}

\begin{proof}
By construction, $W \subseteq V$.  We shall show that if $h \in V$ and $h \perp W$, then $h = 0$.  This will prove $W = V$.

So let $h \in V$ and $h \perp W$.  Then
\begin{equation*}
  h \perp W_0 \plusperp W_1 \plusperp \cdots \plusperp W_m
\end{equation*}
for $m = 0, 1, 2, \cdots$.  By construction,
\begin{equation*}
  V = W_0 \plusperp W_1 \plusperp \cdots \, \plusperp W_m \plusperp V_{m+1}
\end{equation*}
so $h \in V_{m+1}$ for all $m \in \bbz_+$.  By the equations (\ref{eqn5.2}), $h = 0$.
\end{proof}

It will be convenient at this point to introduce some additional standing notation:  Unless otherwise stated, $\bbh = \bbh_0 \supseteq \bbh_1 \supseteq \bbh_2 \supseteq \cdots$ will denote the valuation subspace series for $\bbh$, and $\bbh = H_0 \plusperp H_1 \plusperp H_2 \plusperp \cdots$ will denote the near homogeneous decomposition of $\bbh$.  $V = V_0 \supseteq V_1 \supseteq V_2 \supseteq \cdots$ will denote the valuation subspace series for $V$, and $V = W_0 \plusperp W_1 \plusperp W_2 \plusperp \cdots$ will denote the near homogeneous decomposition of $V$.  $R = R_0 \supseteq R_1 \supseteq R_2 \supseteq \cdots$ will denote the valuation ideal series for $R$.  Thus the meaning of $\bbh_k$, $H_k$, $V_k$, $W_k$, and $R_k$ will be clear without further explanation.

To continue, unless otherwise stated,
\begin{equation}
\begin{split} \label{eqn5.5}
  P_k^\bbh \colon \bbh &\mapsto H_k, \quad Q_{k+1}^\bbh \colon \bbh \mapsto \bbh_{k+1} \\
  P_k^V \colon V &\mapsto W_k, \quad Q_{k+1}^V \colon V \mapsto V_{k+1}
\end{split}
\end{equation}
will denote the orthogonal projections onto $H_k$, $\bbh_{k+1}$, $W_k$, and $V_{k+1}$ respectively. Note that
\begin{equation}
\begin{split} \label{eqn5.6}
  P_k^\bbh + Q_{k+1}^\bbh &= I \text{ on } \bbh_k \\
   P_k^V + Q_{k+1}^V &= I \text{ on } V_k
\end{split}
\end{equation}
where $I$ is the identity map, and
\begin{equation}
\begin{split} \label{eqn5.7}
  \sum_{j=0}^{k} P_j^\bbh + Q_{k+1}^\bbh &= I \text{ on } \bbh \\ 
  \sum_{j=0}^{k} P_j^V + Q_{k+1}^V &= I \text{ on } V
\end{split}
\end{equation}

Now we are equipped to talk about the series of near homogeneous components of individual elements of $\bbh$ and $V$.

\thmcall{definition}\label{def5.8}
Let $h \in \bbh$.  The \emph{near homogeneous decomposition} of $h$ is the convergent series of orthogonal terms
\begin{equation}\label{eqn5.8.1}
  h = P_0^\bbh (h) \plusperp P_1^\bbh (h) \plusperp P_2^\bbh (h) \plusperp \, \cdots
\end{equation}
The individual terms $P_k^\bbh$ are called the \emph{near homogeneous components} of $h$.

If $h \in V$, the \emph{near homogeneous decomposition} of $h$ \emph{with respect} to $V$ is the convergent series of orthogonal terms
\begin{equation}\label{eqn5.8.2}
    h = P_0^V (h) \plusperp P_1^V (h) \plusperp P_2^V (h) \plusperp \, \cdots
\end{equation}
The individual terms $P_k^V$ are called the \emph{near homogeneous components} of $h$ \emph{with respect} to $V$
\exitthmcall{definition}

We shall postpone the proof that the series (\ref{eqn5.8.1}) and (\ref{eqn5.8.2}) are in fact convergent and converge to $h$ until after we establish their connection to the $\ord$ function:

\thmcall{proposition}\label{prop5.9}
For $h \in V$, $\ord(h)$ is the index (subscript) of the first non-vanishing term in the near homogeneous decomposition of $h$ with respect to $V$ if $h \neq 0$, or $\ord(h) = \infty$ if all the terms vanish, that is if $h = 0$.
\exitthmcall{proposition}

\begin{proof}
Suppose $h \in \bbh$ and $h \neq 0$. Then $\ord(h) = m < \infty$, so $h \in V_m = W_m \plusperp V_{m+1}$, and $h \notin V_{m+1}$.  If $P_m^V (h) = 0$, then $h$ would belong to $V_{m+1}$ which is not possible.  Thus $P_m^V (h) \neq 0$.

On the other hand,
\begin{equation*}
  V = W_0 \plusperp W_1 \plusperp \cdots \plusperp W_{m-1} \plusperp V_m
\end{equation*}
so if one or more of the terms $P_j^V (h) \neq 0$ for $j = 0, \dots, m-1$, then $h$ would not belong to $V_m$, which is also not possible.  Thus $P_j^V = 0$, for $j = 0, \dots, m-1$, so the index of the first non-vanishing term in the near homogeneous decomposition of $h$ with respect to $V$ is $\ord(h)$.

Finally, all the terms $P_j^V(h)$ vanish if and only if $h \in \bigintersect V_j$, which is true if and only if $h = 0$ by the equations (\ref{eqn5.2}).  Thus all the terms vanish if and only if $\ord(h) = \infty$.
\end{proof}

\thmcall{corollary}\label{cor5.10}
Let $h \in \bbh$. Then $\ord(h)$ is the index (subscript) of the first non-vanishing term in the near homogeneous decomposition of $h$ if $h \neq 0$, or $\ord(h) = \infty$ if all the terms vanish.
\exitthmcall{corollary}

\begin{proof}
Consider the special case where $\bbh = V$.
\end{proof}

\thmcall{corollary}\label{cor5.11}
If $h$ and $g$ belong to $V$, and $P_k^V (h) = P_k^V (g)$ for all $k \in \bbz_+$, then $h = g$.
\exitthmcall{corollary}

\begin{proof}
$P_k^V (h - g) = 0$ for all $k \in \bbz_+$.  Apply proposition \ref{prop5.9}
\end{proof}

\thmcall{corollary}\label{cor5.12}
Let $h$ and $g$ belong to $\bbh$.  Suppose $P_k^\bbh (h) = P_k^\bbh (g)$ for all $k \in \bbz_+$.  Then $h = g$.
\exitthmcall{corollary}

\begin{proof}
Consider the special case where $\bbh = V$.
\end{proof}

We used the phrase ``\emph{index} of the first non-vanishing term" in the statement of the proposition and its first corollary, even though we could have used the phrase ``\emph{order} of the first non-vanishing term."  The reason is that that if $P_m^V (h) \neq 0$, then $m = \ord(P_m^V (h) )$.  From now on, we shall use the two phrases interchangeably.

One way of interpreting proposition \ref{prop5.9} and corollary \ref{cor5.10} is that all order functions look a lot like the familiar order function on $H^2(E, D^n)$ introduced in section \ref{sec3}, where $\ord(f)$ was the degree of the first non-vanishing homogeneous polynomial in the expansion (\ref{def3.6}) of $f$ as a sum of ordinary homogeneous polynomials.  In fact, in this case, the expansion mentioned in definition \ref{def3.6} is the near homogeneous decomposition of $f$, and this motivates the use of the phrase ``near homogeneous" in the general case of an arbitrary valuation Hilbert module $\bbh$.

Until now, we have been treating the near homogeneous decomposition of $h$ as a formal series.  We are now in a position to prove that it converges to $h$.  At the considerable risk of boring the experts, we do so now:

\thmcall{proposition}\label{prop5.13}
Let $h \in V$.  Then the near homogeneous decomposition (\ref{eqn5.8.2}) of $h$ with respect to $V$ converges in norm to $h$.
\exitthmcall{proposition}

\begin{proof}
By equation (\ref{eqn5.7}),
\begin{equation*}
  h = \sum_{j=0}^k P_j^V (h) \plusperp Q_{k+1}^V (h)
\end{equation*}
By orthogonality, the partial sums
\begin{equation*}
  \norm{ \sum_{j=0}^k P_j^k (h) }^2 = \sum_{j=0}^k \norm{ P_j^V (h) }^2 \leq \norm{h}^2
\end{equation*}
so the partial sums form a Cauchy sequence in $V$.  Because $V$ is closed, the partial sums converge to some element $g \in V$.  Now because $P_k^V$ is a projection and a bounded linear operator on $\bbh$, and because $P_k^V (P_j^V) = 0$ for $j \neq k$,
\begin{equation*}
  P_k^V (g) = P_k^V (\lim_{l \to \infty} \sum_{j=0}^l P_j^V (h) ) = (P_k^V)^2 (h) = P_k^V (h)
\end{equation*}
by orthogonality, so $g = h$ by corollary \ref{cor5.11}.
\end{proof}

\thmcall{corollary}\label{cor5.14}
Let $h \in \bbh$.  Then the near homogeneous decomposition (\ref{eqn5.8.1}) of $h$ converges in norm to $h$.
\exitthmcall{corollary}

\begin{proof}
Consider the special case where $\bbh = V$.
\end{proof}

We turn now to the study of the near homogeneous components $W_0$, $W_1$, $W_2$, $\! \dots \!$ of $V$.  These have a number of striking properties.  First, each is a near homogeneous subspace.  (Think:\@ A linear combination of ordinary homogeneous polynomials of the same degree is a homogeneous polynomial.)  Second, and most importantly for us, the near homogeneous part of an element $h \in W_k$ determines the element $h$ itself.  We shall now make this precise:

\thmcall{definition}\label{def5.15}
A closed linear subspace $W$ of $\bbh$ is \emph{near homogeneous} if the $\ord$ function is constant on $W \comp \{ 0 \}$.
\exitthmcall{definition}

\thmcall{proposition}\label{prop5.16}
The near homogeneous components of $V$ are all near homogeneous subspaces.  In fact, $\ord(h) = m$ for all non-zero $h \in W_m$.
\exitthmcall{proposition}

\begin{proof}
$W_m$ is a closed subspace because it is an orthogonal complement.  That $\ord(h) = m$ follows immediately from proposition \ref{prop5.9}, but here is another short proof: Let $h \in W_m$ and suppose that $\ord(h) \neq m$.  We shall show that $h$ must equal $\{ 0 \}$.  Since $W_m = V_m \orthcomp V_{m+1}$, $h \in V_m$, and so $\ord(h) \geq m$.  Since $\ord(h) \neq m$, $\ord(h)$ must be $> m$, and so $h \in V_{m+1}$.  But $W_m \perp V_{m+1}$, so $h \perp h$.  Thus $h = 0$.
\end{proof}

It can happen that one or more of the near homogeneous components of $V$ may be $\{ 0 \}$.  For example, simply delete a non-zero homogeneous component of a larger subspace to obtain $V$.\\

We now turn to the crucial property that projection of an element $h \in W_k$ onto $H_k$ determines the element $h$ itself.  In fact, this is true for any near homogeneous subspace $W$, not just for the near homogeneous components of $V$:

\thmcall{proposition}\label{prop5.17}
Let $W$ be a near homogeneous subspace of $\bbh$ not equal to $\{ 0 \}$, and let $m$ be the common value of the $\ord$ function on $W \comp \{ 0 \}$.  Let $L_W$ be $P_m^\bbh$ restricted to $W$, so that
\begin{equation}
  L_W \colon W \mapsto P_m^\bbh (W) \subseteq H_m \label{eqn5.17.1}
\end{equation}
Then $L_W$ is an invertible bounded linear transformation.
\exitthmcall{proposition}

\thmcall{proof}
Linearity and boundedness are trivial.  In fact, $L_W = {P_m^\bbh |}_W$ and $\norm{L_W} \leq \norm{P_m^\bbh} = 1$.  To show that $L_W$ is invertible, let $h \in W$ and $h \neq 0$.  Then $m = \ord(h)$, and by corollary \ref{cor5.10}, $P_m^\bbh (h)$ is the first non-vanishing term in the near homogeneous decomposition (\ref{eqn5.8.1}) of $h$.  Thus $L_W (h) = P_m^\bbh (h) \neq 0$, so $\ker(L_W) = \{ 0 \}$, that is $L_W$ is invertible.
\exitthmcall{proof}

In the exceptional case where $W = \{ 0 \}$, the conclusion of the proposition still holds because $h = 0$ is the \emph{only} element of $W$.

An interesting question is, when does $L_W$ have a \emph{bounded} inverse?  By the closed graph theorem, the question is equivalent to, when is $P_m^\bbh (W)$ closed?

The near homogeneous components of $V$ also have an interesting maximality property, which we mention for completeness:

\thmcall{proposition}\label{prop5.18}
The non-zero near homogeneous components of $V$ are all maximal. In fact, if $W_m \neq \{ 0 \}$ and $W'$ is a near homogeneous subspace of $V$ with $W_m \subseteq W'$, then $W_m = W'$. 
\exitthmcall{proposition}

Our first proof of our abstract Beurling's theorem used this proposition, but, as it turns out, we can base a simpler proof on the projection lemmas in the next section.  Hence we shall not need proposition \ref{prop5.18} in this paper,  and so we leave its proof to the reader.

\section{Near Inner Decompositions and the Projection Lemmas.}\label{sec6}
In the one variable case of Beurling's theorem for the closed shift invariant subspaces of $H^2(E, D^1)$, inner functions play a prominent role, and these arise from closed subspaces $W$ with the property that $W \perp z^m W$ for all $m > 0$.  We might call such subspaces \emph{inner} subspaces.

We can easily extend this idea to closed subspaces of a valuation Hilbert module $\bbh$ over a valuation algebra $R$:

\thmcall{definition}\label{def6.1}
A closed subspace $W$ of $\bbh$ is \emph{near inner} if
\begin{equation}
  W \perp R_1 W  \label{eqn6.1.1}
\end{equation}
where $R = R_0 \supseteq R_1 \supseteq R_2 \supseteq \cdots$ is the valuation ideal series of $R$.  An element $h$ of $\bbh$ is \emph{near inner} if the $1$-dimensional subspace generated by it is near inner.
\exitthmcall{definition}

In the $1$-dimensional scalar valued case, a near inner function in $H^2(\bbc, D^1)$ is simply a scalar multiple of an inner function, hence the name ``near inner" in the abstract case.

Our abstract version of Beurling's theorem will involve a sequence of near inner subspaces.  The simple motivating example is

\thmcall{example}\label{example6.2}
Let $V$ be the closed $\poly{\bbc}{z_1, z_2}$ submodule of $H^2(\bbc, D^2)$ generated by the coordinate functions $z_1$ and $z_2$:
\begin{equation}
  V = W_0 \plusperp W_1 \plusperp W_2 \plusperp \cdots  \label{eqn6.2.1}
\end{equation}
where, as usual, the orthogonal direct sum is the near homogeneous decomposition of $V$: $W_0 = \{ 0 \}$, $W_1 = \LH \{z_1, z_2\}$, $W_2 = \LH \{z_1^2, z_1 z_2, z_2^2\}$, etc.  Here, again as usual, $\LH$ denotes the linear hull.
\exitthmcall{example}

Thus we are led to

\thmcall{definition}\label{def6.3}
Let $V$ be a closed subspace of $\bbh$.  A \emph{near inner} decomposition of $V$ is an orthogonal direct sum of subspaces
\begin{equation}
  V = W^0 \plusperp W^1 \plusperp W^2 \plusperp \cdots  \label{eqn6.3.1}
\end{equation}
such that if $k$ and $m \in \bbz_+$, $h \in W^k$, $g \in V$, $r \in R_1$, and $\ord(r h + g) > m$, then
\begin{equation}
  r h + g \perp W^m  \label{eqn6.3.2}
\end{equation}
\exitthmcall{definition}

A few comments are in order.  First, we have denoted the subspaces $W^0$, $W^1$, $W^2,$ $\dots$ using superscripts rather than subscripts because they are not necessarily equal to the near homogeneous components $W_0, W_1, W_2, \, \dots$ of $V$.  Second, each of the subspaces $W^m$ in a near inner decomposition of $V$ is closed because the decomposition is an orthogonal direct sum, and orthogonal complements in closed subspaces are closed.  (We could equally easily have dropped the condition that $V$ be closed, but instead insisted that each subspace $W^m$ be closed.  In other words, a subspace $V$ given by an orthogonal direct sum of closed subspaces is closed.)  Third, if $\ord(h) \geq m$, we may choose $g = 0$ to make $\ord(r h + g) > m$.  Thus $r h \perp W^m$ by equation (\ref{eqn6.3.2}).  In particular, $R_1 W^m \perp W^m$, so each of the closed subspaces $W^m$ is itself near inner.  For future reference, we shall state this last observation as a proposition:

\thmcall{proposition}\label{prop6.4}
Let $W^0 \plusperp W^1 \plusperp W^2 \plusperp \cdot$ be a near inner decomposition of $V$.  Then each of the subspaces $W^m$ is closed and near inner, and if $h \in W^k$, $m \in \bbz_+$, $r \in R_1$, and $\ord(r h) > m$, then $r h \perp W^m$.
\exitthmcall{proposition}

\begin{proof}
See the remarks above.
\end{proof}

It can be hard to verify in practice that a given orthogonal decomposition of $V$ is near inner, but it is often easier to verify that the decomposition satisfies the conclusion of proposition \ref{prop6.4}.  Thus decompositions satisfying this conclusion deserve a name:

\thmcall{definition}\label{def6.5}
Let $V$ be a closed subspace of $\bbh$.  A \emph{weakly near inner} decomposition of $V$ is an orthogonal direct sum of subspaces
\begin{equation}
  V = W^0 \plusperp W^1 \plusperp W^2 \plusperp \cdots  \label{eqn6.5.1}
\end{equation}
such that if $k$ and $m \in \bbz_+$, $h \in W^k$, $r \in R_1$, and $\ord(r h) > m$, then
\begin{equation}
  r h \perp W^m  \label{eqn6.5.2}
\end{equation}
\exitthmcall{definition}

Note that a weakly near inner decomposition is simply one which satisfies the $g = 0$ case of the definition of a near inner decomposition.  Note also that the subspaces in a weakly near inner decomposition are closed (by the remarks above) and are themselves near inner because the definition of a near inner subspace does not involve the auxiliary function $g$.

An interesting question is, when is a weakly near inner decomposition near inner?  In particular, is this always true for near homogeneous decompositions which are weakly near inner?

\

Near inner decompositions are closely related to $R_1$ invariant subspaces.  Recall from definition \ref{def4.9} that $V$ is $R_1$ invariant if $R_1 V \subseteq V$.  Recall also, from the standing notation introduced in section \ref{sec5}, that $V$ is always a \emph{closed} subspace of $\bbh$, that $V = W_0 \plusperp W_1 \plusperp W_2 \plusperp \cdots$ is the near homogeneous decomposition of $V$, and $V = V_0 \supseteq V_1 \supseteq V_2 \supseteq \cdots$ is the valuation subspace series for $V$.

The reader can quickly verify that the near homogeneous decomposition in our simple motivating example \ref{example6.2} is, in fact, a near inner decomposition.  This immediately leads us to the following general proposition:

\thmcall{proposition}\label{prop6.6}
If $V$ is $R_1$ invariant, then the near homogeneous decomposition of $V$ is near inner.
\exitthmcall{proposition}

\begin{proof}
Suppose $V$ is $R_1$ invariant.  Let $h \in W_k$, $g \in V$, and $r \in R_1$.  Then $r h \in V$, so $r h + g \in V$.  If $\ord(r h + g) > m$, then $r h + g \in V_{m+1}$, so $r h + g \perp W_m$.
\end{proof}

The reader might ask, where did we use the hypothesis that $V$ was closed in the above?  The answer is that only closed subspaces are guaranteed to have near homogeneous decompositions.

It would be strikingly convenient if the near inner property were also a sufficient condition for $V$ to be $R_1$ invariant, but the simplest examples show us this is not so:

\thmcall{example}\label{example6.7}
Let $R = \poly{\bbc}{z_1, z_2}$, and let $\bbh = H^2(\bbc, D^2)$, scalar valued $H^2$ of the polydisk $D^2$.  Let the $\ord$ function on $R$ and $\bbh$ be the usual one, the degree of the first non-vanishing term in the power series expansion for $f$.  Let
\begin{equation}
  V = W_0 \plusperp W_1 \plusperp W_2 \plusperp \cdots  \label{eqn6.7.1}
\end{equation}
where
\begin{equation*}
\begin{split}
W_0 = \{ 0 \}, \; W_1 &= \LH \{z_1, z_2\}, \; W_2 = \LH \{z_1^2, z_2^2\}, \\
                  W_3 &= \LH \{z_1^3, z_1^2 z_2, z_1 z_2^2, z_2^3\}, \; \dots
\end{split}
\end{equation*}
 and, as usual, $\LH$ denotes the linear hull.  Then $V$ is closed, the decomposition (\ref{eqn6.7.1}) is the near homogeneous decomposition of $V$, and it is near inner.  However, $V$ is \emph{not} $R_1$ invariant because $W_2$ lacks the needed element $z_1 z_2$.  And, in fact, the closed $R_1$ invariant subspace generated by $V$ is the closed $\poly{\bbc}{z_1, z_2}$ submodule generated by $z_1$ and $z_2$.
\exitthmcall{example}

We need to add an additional property, which we call the \emph{full projection property}, to obtain necessary and sufficient conditions for a closed subspace $V$ to be $R_1$ invariant.  We shall discuss the full projection property in the next section.

To return to the topic of near inner decompositions, the next two \emph{projection lemmas} will provide most of the remaining part of the proof of our abstract Beurling's theorem in section \ref{sec8}.  To set the stage, let
\begin{equation}
  V = W^0 \plusperp W^1 \plusperp W^2 \plusperp \cdots  \label{eqn6.8}
\end{equation}
be a near inner decomposition of $V$.  As before, because the decomposition (\ref{eqn6.8}) is an arbitrary near inner decomposition, and does not necessarily coincide with the near homogeneous decomposition of $V$, we use superscripts rather than subscripts throughout.

The \emph{first projection lemma} is:

\thmcall{lemma}\label{lemma6.9}
Let $r \in R_1$ and $h \in W^k$.  Suppose $f^m = g^0 + g^1 + \cdots g^m$, where $g^0 \in W^0$, $g^1 \in W^1$, $\dots$ $g^m \in W^m$.  Suppose that $\ord(r h - f^m) > m$.  Then the projection of $r h$ on the subspace $W = W^0 \plusperp W^1 \plusperp \cdots \plusperp W^m$ is $f^m$.  Furthermore,
\begin{equation}
  \norm{f^m}^2 = \norm{g^0}^2 + \norm{g^1}^2 + \cdots + \norm{g^m}^2 \leq \norm{r h}^2  \label{eqn6.9.1}
\end{equation}
\exitthmcall{lemma}

\begin{proof} For any element $w \in W$, $w = w^0 + w^1 + \cdots + w^m$, where $w^0 \in W^0$, $w^1 \in W^1$, $\dots$, $w^m \in W^m$.  By the definition of a near inner decomposition, that is definition \ref{def6.3}, $r h - f^m \perp W^j$ for $j = 0, 1, \dots, m$.  Thus
\begin{equation*}
  \inprod{r h - f^m}{w} = \sum_{j=0}^m \inprod{r h - f^m}{w^j} = 0
\end{equation*}
so $r h - f^m \perp W$.  Hence $f^m$ is the projection of $r h$ on $W$.

As for the norm of $f^m$,
\begin{equation*}
\norm{f^m}^2 = \norm{g^0}^2 + \norm{g^1}^2 + \cdots + \norm{g^m}^2
\end{equation*}
because the near inner decomposition (\ref{eqn6.8}) is orthogonal.  Furthermore,
\begin{equation*}
\norm{r h}^2 = \norm{r h - f^m}^2 + \norm{f^m}^2
\end{equation*} because $r h - f^m \perp f^m$.  Thus inequality (\ref{eqn6.9.1}) holds.
\end{proof}

The \emph{second projection lemma} is:

\thmcall{lemma}\label{lemma6.10}
Let $r \in R_1$, and $h \in W_k$.  Consider the formal series
\begin{equation}
  g^0 + g^1 + g^2 + \cdots  \label{eqn6.10.1}
\end{equation}
where $g^0 \in W^0$, $g^1 \in W^1$, $g^2 \in W^2$, $\cdots$, and let $f^m = g^0 + g^1 + \cdots + g^m$.  Suppose that $\ord(r h - f^m) > m$ for $m = 0, 1, 2, \dots$. Then the formal series (\ref{eqn6.10.1}) converges in norm to an element $f$ in $V$, and $r h = f$.
\exitthmcall{lemma}

\begin{proof}
We shall consider the formal series of positive terms
\begin{equation}
  \norm{g^0}^2 + \norm{g^1}^2 + \norm{g^2}^2 + \cdots  \label{eqn6.11}
\end{equation}
By the first projection lemma, the partial sums are
\begin{equation}
  \norm{f^m}^2 = \norm{g^0}^2 + \norm{g^1}^2 + \cdots + \norm{g^m}^2  \leq \norm{r h}^2  \label{eqn6.12}
\end{equation}
Thus the sequence of partial sums $\norm{f^m}^2$, $m = 0, 1, 2, \cdots$ is Cauchy.  Furthermore, if $l < m$, $\norm{f^m - f^l}^2 = \norm{g^{l+1}}^2 + \cdots + \norm{g^m}^2$, so the sequence of partial sums of the formal series (\ref{eqn6.10.1}) is Cauchy in norm.  Because each $f_m \in V$ and $V$ is complete, the formal series (\ref{eqn6.10.1}) converges to an element $f \in V$.

Now let $l \in \bbz_+$.  For all $m > l$, $\ord(r h - f^m) > l$, so $r h - f^m \in \bbh_l$.  By the comments following definition (\ref{def4.3}), $\bbh_l$ is closed, so $r h - f \in \bbh_l$  Thus
\begin{equation*}
r h - f \in \bigintersect_{l=0}^\infty \bbh_l = \{ 0 \}
\end{equation*}
where the right hand equality follows from the equations (\ref{eqn5.2}).  Thus $f = r h$.
\end{proof}

Note that since $r \in R_1$, $\ord(r h) > 0$.  We have included $g^0$ in the above even though $g^0$ must equal $0$ in order for $\ord(r h - f^m)$ to be $> m$.  In fact, $g^0$ through $g^k$ must all be equal to $0$ for the same reason.

\section{The Full Projection Property.}\label{sec7}
In example \ref{example6.7}, we saw that for a closed subspace $V$ to be $R_1$ invariant, it was not sufficient for the near homogeneous decomposition of $V$ to be near inner because some of the components $W_k$ might not have enough elements. The \emph{full projection property}, which we shall define below, guarantees that this can't happen.

Recall from our standing notation introduced in section \ref{sec5} that $\bbh = \bbh_0 \supseteq \bbh_1 \supseteq \bbh_2 \supseteq \cdots$ is the valuation subspace series for $\bbh$, $\bbh = H_0 \plusperp H_1 \plusperp H_2 \plusperp \cdots$ is the near homogeneous decomposition of $\bbh$, and $P_k^\bbh \colon \bbh \mapsto H_k$ is orthogonal projection onto $H_k$.  Also, because we shall be speaking of near homogeneous decompositions, we shall return to using subscripts rather than superscripts.

\thmcall{definition}\label{def7.1}
The near homogeneous decomposition of $V$ has the full projection property if for each $k$ and $m \in \bbz_+$,
\begin{equation}
  P_m^\bbh (r h + g) \in P_m^\bbh (W_m)  \label{eqn7.1.1} ,
\end{equation}
whenever $r \in R_1$, $h \in W_k$, $g \in V$, and $\ord(rh + g) \geq  m$.
\exitthmcall{definition}

\thmcall{proposition}\label{prop7.2}
If $V$ is $R_1$ invariant, then the near homogeneous decomposition of $V$ has the full projection property.
\exitthmcall{proposition}

\begin{proof}
Suppose $V$ is $R_1$ invariant.  Let $k$ and $m \in \bbz_+$, $r \in R_1$, $h \in W_k$, and $g \in V$.  For simplicity, denote $r h + g$ by $f$, and suppose $\ord(f) \geq m$.

Because $V$ is $R_1$ invariant, $r h \in V$ so $f = r h + g \in V$.  Now there are two cases: If $ord(f) > m$, then $P_m^\bbh (f) = 0$ because $f \in \bbh_{m+1}$ and $H_m \perp \bbh_{m+1}$.  Thus $P_m^\bbh \in P_m^\bbh (W_m)$.  On the other hand, if $\ord(f) = m$, then we can write $f = f_m + F_{m+1}$, where $f_m \in W_m$ and $F_{m+1} \in V_{m+1} \subseteq \bbh_{m+1}$.  Since $F_{m+1} \in \bbh_{m+1}$ and $H_m \perp \bbh_{m+1}$, $P_m^\bbh (F_{m+1}) = 0$.  Thus $P_m^\bbh (f) = P_m^\bbh (f_m) \in P_m^\bbh (W_m)$ as required.
\end{proof}

When checking whether the near homogenous decomposition of $V$ has the full projection property, certain simplifications can be made.  First, we never used the fact that $f = r h + g \in V$ in the $\ord(f) > m$ case.  Thus $P_m^\bbh (r f + g) \in P_m^\bbh (W_m)$ always, whenever $\ord(r h + g) > m$, even if $V$ is not $R_1$ invariant.  Second, $P_0^\bbh (r h + g) \in P_0^\bbh (W_0)$ always, even if $V$ is not $R_1$ invariant, because $r \in R_1$, and so $\ord(r h) > 0$. Thus $P_0^\bbh (r h + g) = P_0^\bbh (g) = P_0^\bbh (g_0) \in P_0^\bbh (W_0)$.  Here we have written $g$ as $g_0 + G_1$, where $g_0 \in W_0$ and $G_1 \in V_1 \subseteq \bbh_1$.  Thus $P_0^\bbh (g) = P_0^\bbh (g_0) + P_0^\bbh (G_1) = P_0^\bbh (g_0) + 0 = P_0^\bbh (g_0)$.

We can summarize these remarks in a corollary:

\thmcall{corollary}\label{cor7.3}
Let  $V$ be an arbitrary closed subspace of $\bbh$, possibly not $R_1$ invariant.  Let $k$, $m$, $h$, $r$, and $g$ be as in the proof of the proposition.  Then $P_0^\bbh (r h + g)$ always $\in P_0^\bbh(W_0)$, and $P_m^\bbh (r h + g) \in P_m^\bbh (W_m)$ whenever $\ord(r h + g) > m$.
\exitthmcall{corollary}

\begin{proof}
See the remarks above.
\end{proof}

As in the case of the near inner property, it would be strikingly convenient if the full projection property were also a sufficient condition for $V$ to be $R_1$ invariant.  However, the following simple one-variable example shows that this is not so:

\thmcall{example}\label{example7.4}
Let $R = \poly{\bbc}{z}$, and let $\bbh = H^2(\bbc, D^1)$, scalar valued $H^2$ of the unit disk $D^1$.  Let $a \in D^1$ with $a \neq 0$, and consider the simplest possible non-constant inner function not vanishing at $0$, the simple Blaschke factor
\begin{equation}
  B(z) = \frac{a - z}{1 - \bar a z} \label{eqn7.4.1}
\end{equation}
Let
\begin{equation}
  V = \{ a_0 + a_1 z + a_2 z^2 B(z) + a_3 z^3 B(z) + \cdots \: : \: \sum \abs{a_j}^2 < \infty \} \label{eqn7.4.2}
\end{equation}
Then $V$ is a closed subspace of $H^2(\bbc, D^1)$, the one-dimensional subspaces
\begin{equation}
\begin{split}
  W_0 = \LH \{ 1 \}, W_1 &= \LH \{ z \}, \\
                     W_2 &= \LH \{ z^2 B(z) \}, W_3 = \LH \{ z^2 B(z) \}, \dots  \label{eqn7.4.3}
\end{split}
\end{equation}
are closed, near homogeneous, and mutually orthogonal (because $B$ is inner), and
\begin{equation}
  V = W_0 \plusperp W_1 \plusperp W_2 \plusperp W_3 \plusperp \cdots  \label{eqn7.4.4}
\end{equation}
is the near homogeneous decomposition of $V$.  Clearly, $V$ has the full projection property.  However, $V$ is not $R_1$ invariant (ie.\@ shift invariant).  In fact, the closed shift invariant subspace generated by $V$ is all of $H^2(\bbc, D^1)$.
\exitthmcall{example}

Variations on this example, except for the part about (\ref{eqn7.4.4}) being the near homogeneous decomposition of $V$, are standard in the one-variable theory of $H^2$ spaces.

So the full projection property is not sufficient for $V$ to be $R_1$ invariant.  As our abstract Beurling's theorem, theorem \ref{thm8.1} below, shows, we must add the property that the near homogeneous decomposition of $V$ is near inner to obtain necessary and sufficient conditions for $V$ to be $R_1$ invariant.

\

However, we do have a partial converse to proposition \ref{prop7.2}:

\thmcall{proposition}\label{prop7.5}
Suppose the near homogeneous decomposition of $V$ has the full projection property.  Suppose $k \in \bbz_+$, $r \in R_1$, $h \in W_k$, and $V' = V + \LH \{ r h \}$.  As in definition \ref{def5.3}, let
\begin{equation}
  V' = W_0' \plusperp W_1' \plusperp W_2' \plusperp \cdots  \label{eqn7.5.1}
\end{equation}
be the near homogeneous decomposition of $V'$.  Suppose $m \in \bbz_+$ and $W_m \subseteq W_m'$.  Then $W_m = W_m'$.
\exitthmcall{proposition}

We can summarize the conclusion of proposition \ref{prop7.5} as ``containment implies equality."

Our first proof of our abstract Beurling's theorem used this proposition, but we shall not need it here because we now base our proof on the projection lemmas of section \ref{sec6}.  Thus we shall leave the proof of proposition \ref{prop7.5} to the reader.

\section{The Abstract Beurling's Theorem.}\label{sec8}
Now we can at last state and prove our abstract Beurling's theorem.  It fully characterizes the closed $R_1$ invariant subspaces of valuation Hilbert modules.  For clarity, we shall repeat some of our standing notation in the statement of the theorem and its proof:

\thmcall{theorem}\label{thm8.1}
Let $V$ be a closed subspace of a valuation Hilbert module over a valuation algebra $R$.  Then $V$ is $R_1$ invariant if and only if the near homogeneous decomposition of $V$ is near inner and has the full projection property.
\exitthmcall{theorem}

\begin{proof}
Necessity follows from proposition \ref{prop6.6}, which states that the near homogeneous decomposition of a closed $R_1$ invariant subspace is near inner, and proposition \ref{prop7.2}, which states that a closed $R_1$ invariant subspace has the full projection property.

As for sufficiency, suppose the near homogeneous decomposition of $V$ is near inner and has the full projection property.  Let $r \in R_1$ and $h \in W_k$, where
\begin{equation}
  V = W_0 \plusperp W_1 \plusperp W_2 \plusperp \cdots  \label{eqn8.2}
\end{equation}
is the near homogeneous decomposition of $V$.  We shall show that $r h \in V$.  First, we shall use the full projection property (definition \ref{def7.1}) to show that there are elements
\begin{equation}
  g_0 \in W_0, \: g_1 \in W_1, \: g_2 \in W_2, \: \dots  \label{eqn8.3}
\end{equation}
such that for $m = 0, 1, 2, \dots$,
\begin{equation*}
  \ord(r h - f_m) > m
\end{equation*}
Here, we have denoted the sum $g_0 + \cdots g_m$ by $f_m$.  Then we shall use the second projection lemma (\ref{lemma6.10}) to conclude that the formal sum
\begin{equation}
  g_0 + g_1 + g_2 + \cdots  \label{eqn8.4}
\end{equation}
converges to an element $f \in V$, and that $r h = f$.  This will prove that $r h \in V$.

First, let us construct the sequence (\ref{eqn8.3}) by induction.  Because $r \in R_1$, $\ord(r h) > 0$.  Thus $P_0^\bbh (r h) = 0$, so we can (and must) choose $g_0 = 0$.  Now suppose $g_0, \dots, g_{m-1}$ have been constructed inductively, so that $g_j \in W_j$ and $\ord(r h - f_j) > j$ for $j = 0, \dots, m - 1$.  We shall use the full projection property of the near homogeneous decomposition (\ref{eqn8.2}) of $V$ to construct $g_m$.  Note that $\ord(r h - f_{m-1}) \geq m$.  By the full projection property, $P_m^\bbh (r h - f_{m-1}) \in P_m^\bbh (W_m)$.  If $W_m \neq \{ 0 \}$, we may use proposition \ref{prop5.17} to conclude that there is a unique element $g_m \in W_m$ such that $P_m^\bbh (g_m) = P_m^\bbh (r h - f_{m-1})$.  In the exceptional case where $W_m = \{ 0 \}$, we can (and must) choose $g_m = 0$.  (See the remark following proposition \ref{prop5.17}.)  Let $f_m = f_{m-1} + g_m$.  Then $f_m = g_1 + \cdots + g_m$, and $g_m \in W_m$ by construction.

But we still must prove that $\ord(r h - f_m) > m$.  To that end, recall from corollary \ref{cor5.10} that $\ord(r h - f_m)$ is the index of the first non-vanishing term in the near homogeneous decomposition of $r h - f_m$, or $\infty$ if $r h - f_m = 0$.  Since $\ord(r h - f_{m-1}) > m - 1$, all the terms in the near homogeneous decomposition of $r h - f_{m-1}$ vanish up to and including the index $m - 1$, that is $P_j^\bbh (r h - f_{m-1}) = 0$ for $j = 0, \dots, m - 1$.  Also, because $g_m \in W_m$, $\ord(g_m) = m$ (or $\infty$ if $g_m = 0$).  Thus $P_j^\bbh (g_m) = 0$ as well for $j = 0, \dots m-1$.  Furthermore, $P_m^\bbh (r h - f_{m-1}) = P_m^\bbh (g_m)$ by construction.  Thus
\begin{equation*}
P_j^\bbh (r h - f_m) =  P_j^\bbh (r h - f_{m-1}) - P_j^\bbh (g_m) =0
\end{equation*}
for $j = 0, \dots, m$, and so $\ord(r h - f_m) > m$.

This concludes the first part of the proof that $r h \in V$.  The second part is immediate.  Because the near homogeneous decomposition (\ref{eqn8.2}) of $V$ is near inner, we may apply the second projection lemma.  By the second projection lemma, the formal series \ref{eqn8.4} converges in norm to an element $f \in V$, and $r h = f$.  Thus $r h \in V$.

Now let $h \in V$, not necessarily an element of $W_k$.  Let $r \in R_1$ as before, and let
\begin{equation}
  h = h_0 \plusperp h_1 \plusperp h_2 \plusperp \cdots  \label{eqn8.5}
\end{equation}
be the near homogeneous decomposition of $h$ with respect to $V$.  By definition \ref{def5.8}, each $h_k \in W_k$, and by proposition \ref{prop5.13}, the series (\ref{eqn8.5}) converges in norm to $h$.  By the definition of a Hilbert module, the linear operators $m_r \colon f \to r f$, $f \in \bbh$ are bounded in operator norm, so the series
\begin{equation*}
  r h = r h_0 + r h_1 + r h_2 + \cdots
\end{equation*}
converges in norm to $r h$.  Further, each term $r h_k \in V$, so, since $V$ is closed, $r h \in V$ as required.  Thus $V$ is $R_1$ invariant.
\end{proof}

The last part of this proof is useful enough to isolate as a separate corollary:

\thmcall{corollary}\label{cor8.6}
A closed subspace $V$ is $R_1$ invariant if and only if, for each $k \in \bbz_+$ and $r \in R_1$, $r W_k \subseteq V$.
\exitthmcall{corollary}

\begin{proof}
See the above.
\end{proof}

A proper generalization of Beurling's theorem should characterize, not just closed $R_1$ invariant subspaces of $\bbh$, but closed $R$ submodules as well.  This we can do in two important cases.  Recall from definition \ref{def4.10} that a valuation algebra $R$ is \emph{unital} if $R$ has a unit element $1_R$ and
\begin{equation}
  R = \{ \lambda 1_R : \lambda \in \bbc \} + R_1  \label{eqn8.7}
\end{equation}
From proposition \ref{prop4.11} we immediately have

\thmcall{corollary}\label{cor8.8}
Let $R$ be a valuation algebra satisfying either of the following conditions:
\enumcall{enumerate}
  \item $R = R_1$, or
  \item $R$ is unital
\exitenumcall{enumerate}
Then a closed subspace $V$ of a valuation Hilbert module $\bbh$ over $R$ is an $R$ submodule of $\bbh$ if and only if the near homogeneous decomposition of $V$ is near inner and has the full projection property.
\exitthmcall{corollary}

\begin{proof}
Apply proposition \ref{prop4.11} and theorem \ref{thm8.1}.
\end{proof}

We remind the reader that in the iconic case where $R = \poly{\bbc}{z_1, \dots, z_n}$, $R$ is unital, so corollary \ref{cor8.8} fully characterizes the closed $\poly{\bbc}{z_1, \dots, z_n}$ submodules of $\bbh$.

\

\section{Example -- $H^2$ Spaces.}\label{sec9}
The closed invariant subspaces of many familiar Hilbert spaces of analytic functions are completely described by our Abstract Beurling's Theorem, Theorem\ref{thm8.1} and its corollaries. In this section, we shall describe the closed invariant subspaces of $\bbh = H^2(E, D^n)$, where, as usual, $E$ is a complex Hilbert space and $D^n$ is the $n$-dimensional polydisk.  Again, as usual, $R$ will be the familiar algebra of complex valued polynomials in the variables $z_1, \ldots, z_n$.  And both $R$ and $\bbh$ will have the usual ord function, $\ord(f) =$ the multiplicity of the zero of f at $0$.

\thmcall{example}\label{example9.1}
Theorem \ref{thm8.1} and its corollaries completely characterize the closed invariant subspaces of $\bbh = H^2(E, D^n$).  
\exitthmcall{example}

\begin{proof}
All we have to show is that the $\ord$ function on $\bbh$ is upper semi-continuous. So let $f_j$, $j = 1, 2, 3, \dots$ be a norm convergent sequence in $\bbh$ converging to $f$.  For what follows, we need to briefly review one more thing about multi-index notation:  Recall  that $\abs{m} = \abs{m_1} + \cdots \abs{m_n}$.

Now by definition \ref{def3.5}, the lowest degree terms in the power series expansion of $f$ have degree $\ord(f)$. In multi-index notation, let $a_m z^m$ be one of these terms (so $\ord(f) = \abs{m}$), and let $a_{j, m}$be the corresponding term in the power series expansion of $f_j$.  Then

\begin{equation*}
   a_m = \int_0^{2\pi} f(e^{\imath \theta}) e^{-\imath m \theta} \dtheta/(2 \pi^n)
\end{equation*}
in multi-index notation, with a corresponding formula for $a_{j, m}$ with $f$ replaced by $f_j$.  Here, the integral is repeated $n$ times, $\theta = \theta_1 \cdots \theta_n$, and $\dtheta = \dtheta_1 \cdots \dtheta_n $.

The reader will note that $a_m$ is simply the $m^\text{th}$ Fourier coefficient of $f$.  Since the sequence $f_j$ converges to $f$ in $H^2$ norm (and $e^{-\imath m \theta}$ is bounded)  $a_{j, m}$ converges to $a_m$ as $j \rightarrow \infty$. Thus, for $j \geq \text{ some } J$, $\abs{a_{j, m}} \geq \abs{a_m}/2 \neq 0$. Hence, for each $j \geq J$, $f_j$ has terms in its power series expansion of degree $\abs{m} = \ord(f)$, so for $j \geq J$, $\ord(f_j) \leq \ord(f)$. Consequently, $\limsup \ord(f_j) \leq \ord(f)$, so the $\ord$ function is upper semi-continuous on $\bbh$.
\end{proof}

In a second paper, we shall use theorem \ref{thm8.1} to fully characterize the closed invariant subspaces of a number of Hilbert spaces of analytic functions in several complex variables. These will include $H^2$ of the unit ball in $\bbc^n$, and of general bounded symmetric domains. Also, these will include weighted Bergman spaces of analytic functions on domains in $\bbc^n$, and weighted Bergman spaces of analytic functions and analytic differentials on complex analytic manifolds.

\end{document}